\documentclass[12pt]{amsart}

\usepackage{fullpage}
\usepackage{amsmath}
\usepackage{amssymb}
\usepackage{color}
\usepackage{epsfig}
\usepackage{subfig}

\newtheorem{theorem}{Theorem}[section]

\newtheorem{definition}[theorem]{Definition}
\newtheorem{lemma}[theorem]{Lemma}

\newcommand{\emptyroom}{\varepsilon}
\newcommand{\mcP}{\mathcal{P}}
\newcommand{\bbC}{\mathbb{C}}
\newcommand{\bbZ}{\mathbb{Z}}

\usepackage{setspace}
\begin{document}

\begin{abstract} 
We verify a recent conjecture of Kenyon/Szendr\H{o}i by computing the generating function
for pyramid partitions.  Pyramid partitions are closely related to
Aztec Diamonds; their generating
function turns out to be the partition function for the Donaldson-Thomas
theory of a non-commutative resolution of the conifold singularity $\{x_1x_2 -
x_3x_4 = 0 \} \subset \bbC^4$.  The proof does not require algebraic geometry; it uses a modified version of the domino
shuffling algorithm of Elkies, Kuperberg, Larsen and Propp, ~\cite{EKLP_II}.
\end{abstract}

\title{Computing a Pyramid Partition generating function with dimer shuffling}
\author{Ben Young}
\maketitle

\section{Introduction}

Consider the pyramid-shaped stack of square bricks shown in Figure~\ref{fig:bricks_and_room}.  The bricks are the same ones used to $q$-enumerate Aztec Diamonds in~\cite{EKLP_II}: ridges on the top and bottom of the bricks restrict the manner in which the bricks may be stacked.  Each brick rests upon two side-by-side bricks, and is rotated 90 degrees from the bricks immediately below it.  We use two colors of bricks --- light and dark --- to make alternating layers of this pyramid, starting with dark bricks at the pyramid's apex.

In Figure~\ref{fig:bricks_and_room}, there is a row of three dark bricks at the top of the pyramid.  It is straightforward to build a similar pyramid with a row of $n \geq 1$ bricks along the top.  Following~\cite{NC_DT}, we make the following definitions:

\begin{definition}
\label{defn:emptyroom} The pyramid with a row of $n$ dark bricks at the top is called \emph{the empty room\footnote{This admittedly strange terminology is borrowed from the jargon of 3D partitions, which are made of stacks of boxes in the corner of a room.  Here, the configuration of minimum weight is an empty room, with no boxes.} of length $n$}, and is denoted $\emptyroom_n$.
\end{definition}

\begin{definition}
\label{defn:pyramid_partition}
A \emph{pyramid partition of length $n$} is a finite subset $\pi$ of the bricks of $\emptyroom_n$ such that if $B$ is a brick in $\pi$, then all of the bricks of $\emptyroom_n$ which rest upon $B$ are also in $\pi$.  Let $\mcP_n$ denote the set of all pyramid partitions of length $n$.
\end{definition}

\begin{definition}
\label{defn:partition_weight}
The \emph{weight} of $\pi$, $w_0(\pi)$, is
\[
q_0^{\#\{\text{dark bricks in }\pi\}}
q_1^{\#\{\text{light bricks in }\pi\}}.
\]
\end{definition}

\begin{figure}
\caption{A pyramid partition of length 1, viewed from the side and from above}
\label{fig:two_views}
\input{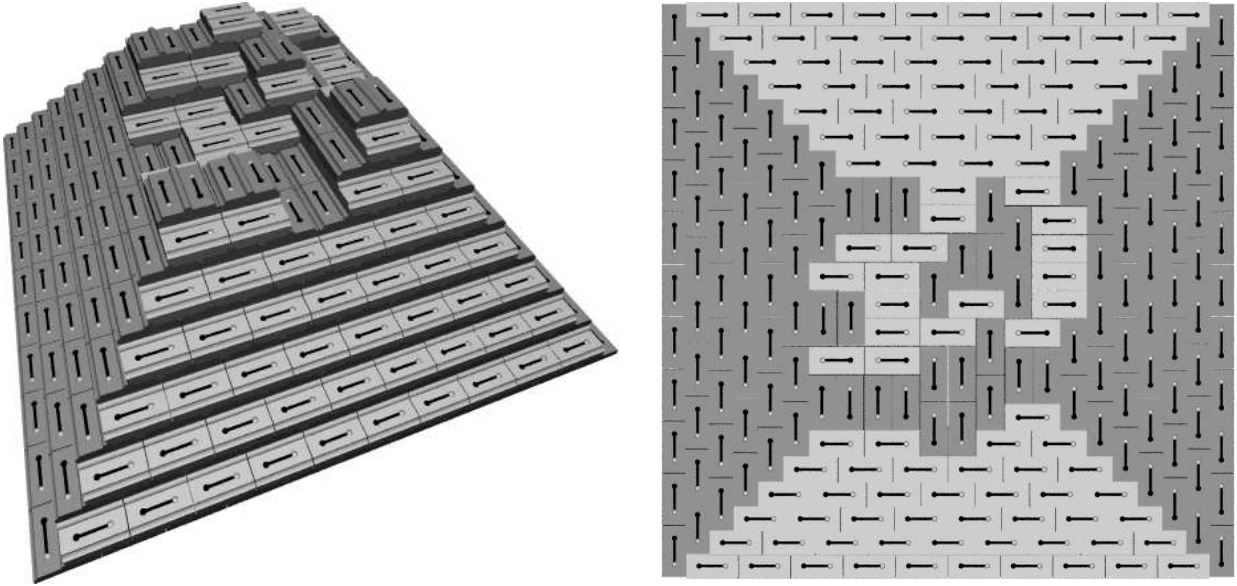}
\end{figure}

In other words, a pyramid partition is a collection of bricks removed from
$\emptyroom_n$ such that the remaining pile of bricks is stable.  For our
treatment, it is better to draw pyramid partitions by drawing the remaining
pile of bricks.  For an example of a pyramid partition drawn in this way, see
Figure~\ref{fig:two_views}.   Note that $\emptyroom_n$ is itself a pyramid partition of weight 1, for all $n$.

There is a third way to view a pyramid partition $\pi$, which is much more
useful computationally.  Recall that a \emph{dimer cover} (or \emph{1-factor})
of a graph $G$ is a subgraph $G'$ such that every vertex of $G'$ has degree 1.
Each brick in $\pi$ has two dimers stencilled on the top; dark bricks have vertical (North-South) dimers, whereas light bricks have horizontal (East-West) dimers.  When one views $\pi$
from above, one can see a dimer cover of the square lattice (see
the right-hand image in Figure~\ref{fig:two_views}).  It is helpful to think of
the lattice points as pairs of half-integers, so that the origin lies above the
axis of symmetry of $\emptyroom_n$.  

\begin{figure}
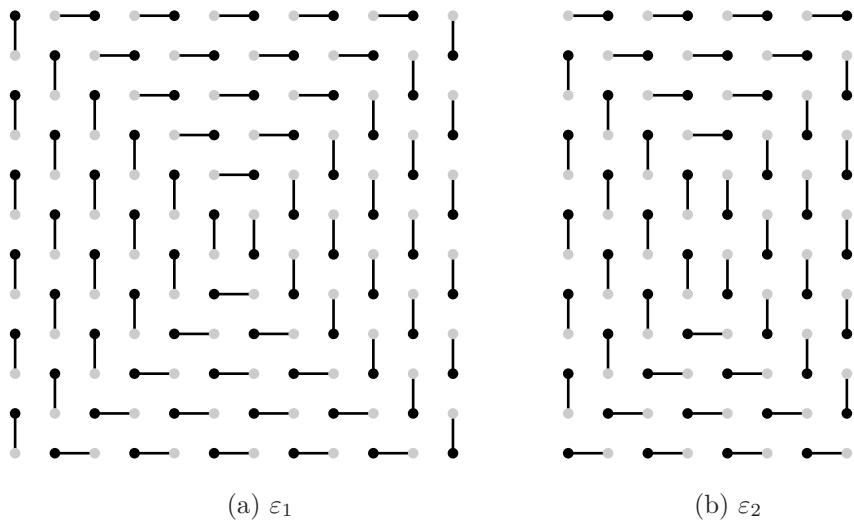

\caption{The empty rooms of lengths 1 and 2.}
\label{fig:empty_room_1_2}
\subfloat[$\emptyroom_1$]{\input{empty_room_1.pstex_t} }
\subfloat[$\emptyroom_2$]{\input{empty_room_2.pstex_t} }
\end{figure}

Since every pyramid partition has only
finitely many bricks, the dimer cover associated to $\pi$ looks like that of $\emptyroom_n$ (see Figure~\ref{fig:empty_room_1_2}) when one moves far enough from the origin.  Indeed, given a dimer cover $T$ of the square lattice which is asymptotically identical to $\emptyroom_n$, it is straightforward to construct a corresponding pyramid partition which looks like $T$ from above.  We shall therefore refer to these dimer configurations as pyramid partitions, as well.

\begin{figure}
\caption{Special bricks, assembled into the configuration $\emptyroom_3$. }
\label{fig:bricks_and_room}
\input{bricks_and_room.pstex_t}
\end{figure}

In~\cite{NC_DT}, Szendr\H{o}i defines a bivariate generating function for $\mcP_n$ by
\[
Z_A^{(n)}(q_0,-q_1) = \sum_{\pi \in \mcP_n}w_0(\pi)
\]
and observes that $Z_A^{(1)}(q_0, q_1)$ arises as the partition function for the Donaldson-Thomas theory of a non-commutative resolution of the conifold singularity $\{x_1x_2 - x_3x_4 = 0 \} \subset \bbC^4$.  
Szendr\H{o}i conjectures that 
\begin{equation}
\label{eqn:balazs_conjecture}
Z_A^{(n)}(q_0,-q_1) = M(1, q_0q_1)^2 
\prod_{k \geq 1}(1+q_0^kq_1^{k-1})^{k+n-1}
\prod_{k \geq 1}(1+q_0^kq_1^{k+1})^{\max(k-n+1,0)}
\end{equation}
where $M(x,q)$ is the MacMahon function
\[M(x,q) = \prod_{n=1}^{\infty}\left(\frac{1}{1-xq^n}\right)^n.\]
This conjecture (or at least the special case $q_0=q_1=q$) was originally posed by Kenyon~\cite{KENYON} 

We present a proof of this conjecture. We first do the case $n=1$, using a modification of the domino shuffling argument of \cite{EKLP_II}, originally used to compute the weight generating function of an Aztec Diamond.  Strikingly, this case uses the Donaldson-Thomas partition function of the resolution of this conifold, computed in~\cite{SR_DT}.

Before we go any further, let us choose a more convenient notation.
\begin{definition}
Let 
\begin{align*}
Z(n; q_0,q_1) & := Z^{(n)}_A(q_0, -q_1) = \sum_{\pi \in \mcP_n}w_0(\pi); \\
Z(\infty; q_0, q_1)& := M(1, q_0q_1)^2M(-q_1^{-1}, q_0q_1)^{-1}.
\end{align*}
\end{definition}

We may now restate (and prove) Equation (\ref{eqn:balazs_conjecture}) for $n=1$ in the following form:
\begin{theorem}
\label{thm:main}
$Z(1; q_0,q_1) = M(-q_1^{-1}, q_0q_1)^{-1} Z(\infty; q_0, q_1).$
\end{theorem}

We have chosen the notation somewhat suggestively here.  Our proof, \emph{very} informally speaking, is that domino shuffling transforms pyramid partitions of length $n$ into pyramid partitions of length $n+1$ in a weight-preserving manner (the transformation is not quite bijective). Repeating this procedure forever, we get ``pyramid partitions of length $\infty$''.  These objects are easily weight-enumerated due to a surprising bijection with a type of $3D$ partitions which we have called \emph{super-rigid partitions} (see Section~\ref{sec:super_rigid}). 
It is also possible to use our methods to prove equation~(\ref{eqn:balazs_conjecture}) for general $n$, which in our new notation looks like this:
\begin{equation}
\label{thm:generalization_of_main}
Z(n; q_0,q_1) = M(1, q_0q_1)^2 
\prod_{k \geq 1}(1+q_0^kq_1^{k-1})^{k+n-1}
\prod_{k \geq 1}(1+q_0^kq_1^{k+1})^{\max(k-n+1,0)}
\end{equation}
In section 7, we shall outline how to modify our proof of Theorem~\ref{thm:main} to handle this more general case.  The proof is relegated to a later section of the paper because it contains essentially no new combinatorial ideas (only greater complication) and because the $n=1$ case is of greater geometric interest.

\section{Dimer Shuffling}

Next, we will describe the shuffling algorithm, originally published in
\cite{EKLP_II}.  We shall call this algorithm \emph{dimer shuffling}, rather than domino shuffling, since all of our pictures are of dimers, which are dual to the dominos of \cite{EKLP_II}.  However, the shuffling algorithm is identical.  We review it here in order to define all of our terminology.

The purpose of the algorithm is to transform a pyramid partition of length $n$
into a pyramid partition of length $n+1$.  Unfortunately, the dimer shuffle is
not quite an honest function from $\mcP_n$ to $\mcP_{n+1}$, in that there are
several different possible outcomes of the algorithm.  So let us first describe the deterministic part of the algorithm, the \emph{sliding map}, which acts on certain partial dimer covers $T$ of the square lattice.

First of all, we colour the vertices of the lattice black and white in a
checkerboard pattern.  Any dimer on this lattice has one endpoint of each
color.  Of course, we must pick the parity of this colouring; it depends on the
parity of $n$ (see Figure~\ref{fig:empty_room_1_2}).  If $n$ is odd, then the
center square of the lattice has a black vertex in the upper left corner.
Otherwise, that vertex is white. 

We adopt the following definitions of \cite{EKLP_II} (changing the notation slightly):
\begin{definition}
Two side-by-side dimers (or, sometimes, their four endpoints) are called a \emph{block}.   A block is \emph{odd} if it has a black vertex in the upper left corner; otherwise it is \emph{even}.  
\end{definition}

Figure~\ref{fig:odd_even}(a) shows the different types of odd and even blocks. 
As you can see in Figure~\ref{fig:empty_room_1_2}, the empty room of length $n$ always has precisely $n$ odd blocks in a vertical line in the center.

\begin{figure}
\caption{}
\label{fig:odd_even}
\rule{0in}{\baselineskip}
\subfloat[Odd blocks and even blocks]{ \input{odd_even.pstex_t} }
\subfloat[The directions in which dimers move during sliding] { \input{dimer_directions.pstex_t} }
\end{figure}

\begin{definition}
\label{def:odd_deficient}
An \emph{odd-deficient} (respectively, \emph{even-deficient}) dimer cover is a partial dimer cover such that the set of non-covered vertices is a finite union of odd (respectively, even) blocks.  
Given a dimer cover $T$,  construct the odd-deficient dimer cover $\tilde T$ by deleting all of the odd blocks of $T$.  Construct the even-deficient dimer cover $\hat T$ by deleting all of the even blocks of $T$.  Let
\begin{align*}
\tilde{\mcP}_n &:= \{ \tilde{\pi} : \pi \in \mcP_n \}, \\
\hat{\mcP}_n &:= \{ \hat{\pi} : \pi \in \mcP_n \}.
\end{align*}
\end{definition}

\begin{definition}
\label{def:sliding_map}
The \emph{sliding map} $S$ is a mapping from the set $\{$dimers on the colored square lattice$\}$ to itself.  If $d$ is a dimer, then define $S(d)$ to be the other dimer in the odd block containing $d$.
If $T$ is an odd-deficient partial dimer cover, then define $S(T)$ to be the partial dimer cover $\{S(d)\;:\;d \in T\}$.
\end{definition}

Observe that $S$ moves each dimer in $T$ one unit to the north, south,
east, or west, depending on its position; Figure ~\ref{fig:odd_even}(b) shows
the directions in which the dimers move.  We shall often call dimers
\emph{northbound, southbound, eastbound, or westbound}, according to the
direction in which they slide.  
Note that $S$ depends on the parity of the lattice coloring we have chosen.

\begin{lemma}
$S$ is an involution on the set of odd-deficient dimer covers. The restriction $S|_{\tilde \mcP_n}$ is a bijection from $\tilde \mcP_n$ to $\hat \mcP_{n+1}$ with their usual colorings.
\end{lemma}

\begin{proof}
One first shows that $S$ is an involution, essentially by analyzing all of the possible local odd-deficient configurations of dimers.  This is done in detail in~\cite{EKLP_II}.  To verify that the image of $S$ is $\hat \mcP_{n+1}$, observe that $S(\tilde{\emptyroom}_n) = \emptyroom_{n+1}.$   The parity of the usual coloring of $\mcP_{n+1}$ is the opposite of that of $\mcP_n$, so for $\pi \in \mcP_n$, $S(\pi)$ is even-deficient and asymptotic to $\emptyroom_{n+1}$. 
\end{proof}

Figure~\ref{fig:shuffling_works} shows how $S$ works. In (a), we have deleted
all of the odd blocks of the pyramid partition in Figure~\ref{fig:two_views};
the missing odd blocks are marked with grey squares.  In (b), we have applied
$S$, and now the grey squares denote the missing \emph{even} blocks.
Observe that $S(\pi) \in \hat \mcP_2$.

\begin{figure}
\caption{}
\label{fig:shuffling_works}
\subfloat[An odd-deficient $\tilde{\pi} \in \tilde{\mcP}_1$.]{ \input{shuffle_example.pstex_t} }
\subfloat[$S(\tilde{\pi})$.]{ \input{shuffle_example2.pstex_t} }
\end{figure}


We may now define the dimer shuffling algorithm, which extends $S$ to a map
\[S:\mcP_n \rightarrow \{ \text{formal sums of pyramid partitions of length } n+1\}.\]

\begin{definition} 
Let $\pi \in \mcP_{n}$.  The following three steps constitute the \emph{dimer shuffling algorithm}:
\begin{enumerate}
\item(Deleting) Delete all of the odd blocks in $\pi$ to get $\tilde{\pi}$. 
\item(Sliding) Compute $S(\tilde{\pi})$, as defined above.
\item(Creating) Now we have a partial dimer cover which is possibly missing some \emph{even} blocks.  Each block may be filled in with either two horizontal dimers, or two vertical dimers.  Define $S(\pi)$ to be the formal sum of all of these fillings.
\end{enumerate}
\end{definition}

It is fairly straightforward to see that these steps are well-defined and that they do indeed give you a formal sum of dimer covers of the plane; this is shown in detail in~\cite{EKLP_II}.

Finally, let us prove a lemma about the number of odd blocks of a pyramid partition.  Observe that Figure~\ref{fig:shuffling_works}(a) has 10 odd blocks, whereas (b) has 9 even blocks.  In general, we have:

\begin{lemma}
\label{lem:odd_block_count}
Let $\tilde \pi \in \tilde \mcP_n$. Then
$ \#\{\text{odd blocks in } \tilde \pi \} - \#\{\text{even blocks in } S(\tilde \pi) \} = n.$
\end{lemma}

\begin{proof} Suppose there are $m$ odd blocks in $\pi$ and $m'$ even blocks in $S(\pi)$.  Let $R$
be a $2a \times (2a+2n-2)$ rectangle of lattice points centered at the origin,
where $a$ is large enough that $\pi$ is identical to $\emptyroom_n$ outside
$R$, and there are no odd blocks of $\pi$ on the boundary of $R$.  For
example, for the odd-deficient partition of
Figure~\ref{fig:shuffling_works}(a), we could take $a=7$ and $R$ to be the $14
\times 14$ rectangle of lattice points shown in the illustration.   

Each dimer has two endpoints and each (missing) odd block has four vertices, so the number of dimers in $R$ is
\begin{equation}
\label{eqn:block_count_1}
\frac{(2a)(2a+2n-2) - 4m}{2}.
\end{equation}
Now let us shuffle the dimers in $R$.  The same dimers now fit into a $(2a -2) \times (2a + 2n)$ rectangle, which has $(2a - 2)(2a+2n)$ and contains all $m'$ odd blocks.  So the number of dimers in $R$ is also equal to
\begin{equation}
\label{eqn:block_count_2}
\frac{(2a-2)(2a+2n) - 4m'}{2}
\end{equation}
Setting Equations~(\ref{eqn:block_count_1}) and~(\ref{eqn:block_count_2}) equal, we obtain the lemma. 
\end{proof}

\section{Weighting the lattice}

In order to use domino shuffling as a computational tool, we need to find a way to calculate the weight of a pyramid partition from its dimer form, without interpreting it as a pile of bricks.  Our strategy shall be to assign a monomial weight to every edge of the square lattice  
 in such a way that the renormalized product of the edge weights of any pyramid partition $\pi$ is $w_0(\pi)$.  This idea is mentioned in~\cite{NC_DT},  but we shall need to be explicit about what edge weights we use and how we do the renormalization.

In order to determine the proper weights to use, it is helpful to consider how a minimal change in the dimer configuration should affect the weight.  We make the following definition:

\begin{definition}
Let $\pi$ be a pyramid partition.  An \emph{elementary move} is the act of adding an appropriately colored block to $\pi$ to obtain a new pyramid partition.
\end{definition}

\begin{figure}
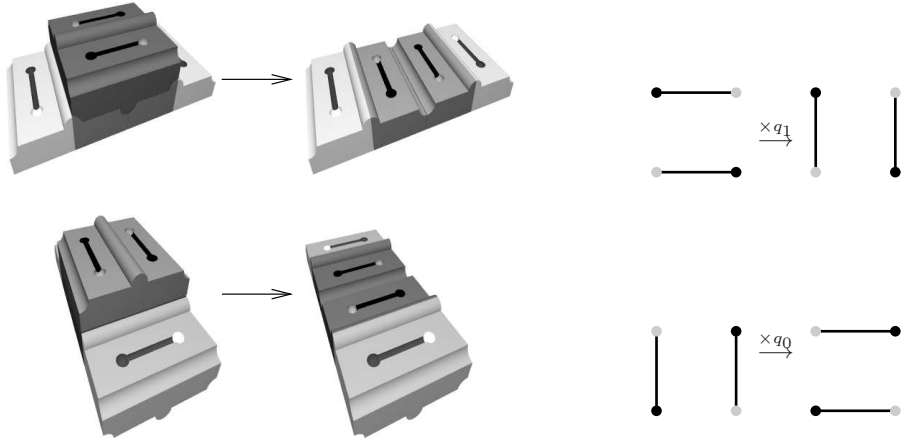

\caption{Elementary moves for adding a bricks to pyramid partitions}
\label{fig:elementary_move}
\input{elem_brick_moves.pstex_t}
\rule{0.5in}{0in}
\input{elementary_move.pstex_t}
\end{figure}

When we analyze the effect an elementary move has on the dimer version of $\pi$, we see that there are two different types of elementary moves for adding a dark or light brick.   They are shown in Figure~\ref{fig:elementary_move}; recall that our convention in drawing the brick pictures is to show the \emph{complement} of the pyramid partition!  An odd elementary move should contribute $q_0$ to the weight, whereas an even move should contribute $q_1$.

We may now assign a weight to each edge of the square lattice which is compatible with the elementary moves, in the following sense: select any $2\times2$ block of vertices in the weighted lattice.  If it is an odd block, we should have
\[\frac{\text{weight of two horizontal dimers}}{\text{weight of two vertical dimers}} = q_0,\]
and if it is an even block, we should have
\[\frac{\text{weight of two vertical dimers}}{\text{weight of two horizontal dimers}} = q_1.\]

In fact, there are many ways to do this, but it is convenient to
choose the weighting in which all vertical edges have weight 1, and all the
northbound horizontal edges closest to the $x$ axis have weight 1 (see
Figure~\ref{fig:weighted_lattice}).  We adopt the convention that in a weighted
lattice, edges with no marked weight get weight 1.

\begin{definition}
If $d$ is a dimer, then $w_0(d)$ is the weight assigned to $d$ in Figure~\ref{fig:weighted_lattice}.
\end{definition}

Now we need to explain how to use these edge weights to compute the weight
of a pyramid partition $\pi$.  Naively, we want to say that the weight of
$\pi$ is the product of the weights of its edges.  However, since $\pi$ covers
the entire plane and has an infinite number of edges, this is meaningless.  Fortunately, all one has to do is to normalize the weight in the following sense:

\begin{figure}
\caption{The $w_0$ weighting on the square lattice.  The heavy black line is the $x$ axis.}
\label{fig:weighted_lattice}
\rule{0in}{\baselineskip} \\
{
\tiny
\input{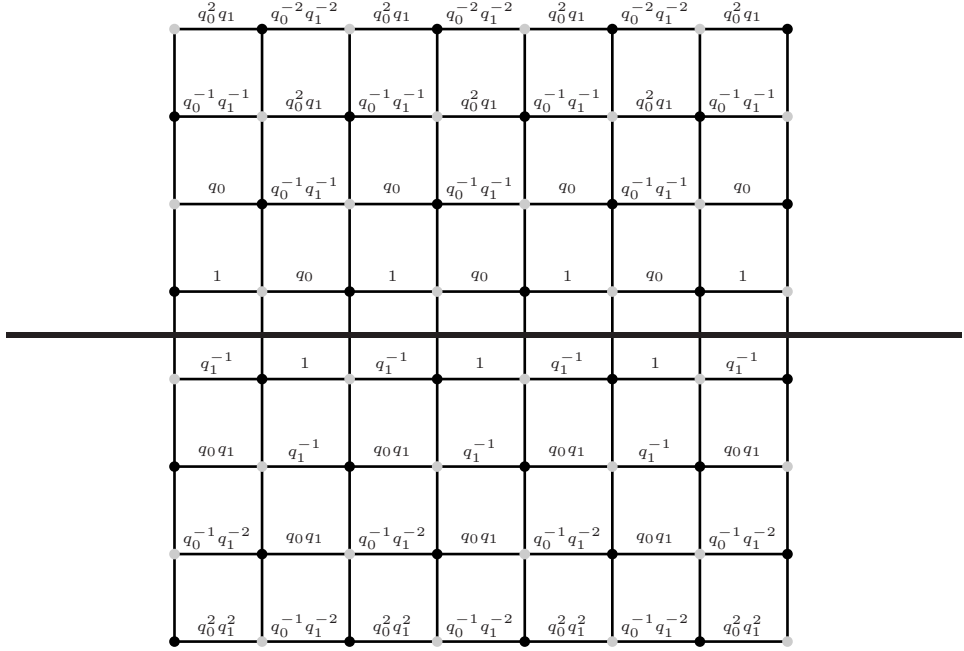}
}
\end{figure}

\begin{lemma} Suppose that $\pi \in \mcP_n$.  Let $R$ be a finite region of the lattice which contains all of the edges where $\pi$ differs from $\emptyroom_n$. Then 
\[
w_0(\pi) = 
\left(
\prod_{e \in R \cap \pi} w_0(e)
\right) \left(
\prod_{e \in R \cap \emptyroom_n} w_0(e)
\right)^{-1}.
\]
\end{lemma}

\begin{proof}
As a base case, let $\pi = \emptyroom_n$ and observe that both sides are equal to 1.  Next, suppose that the lemma holds for some pyramid partition $\pi_0$; by the preceding remarks, it also holds for all $\pi$ which differ from $\pi_0$ by an elementary move.  The lemma then follows by induction on the number of bricks in $\pi.$
\end{proof}

\section{Weighting and Shuffling}

We shall use a different weighting function, $w_1$, to weight $S(\pi)$.  Essentially, we want to think of the weight of a dimer as being unaffected by the shuffling operation.  In fact, we shall define a series of weight functions $w_1, w_2, w_3, \ldots$, which have the property that \[w_0(d) = w_1(S(d)) = w_2(S^2(d)) = \cdots\] for any dimer $d$.
\begin{definition}  
Let $d$ be dimer in a pyramid partition of length $n$ (with the usual lattice coloring).  Let $a \geq 1$.  Define the weight function $w_a$ by
\[
w_a(d) = w_0(\underbrace{S^{-1} \circ S^{-1} \circ \cdots \circ S^{-1}}_{a-1}(d)).
\]
\end{definition} 

For a comparison of $w_0$ and $w_1$, see Figure~\ref{fig:weighted_lattice}.  Observe that if $d$ is a vertical dimer, then $w_a(d) = 1$ for all $a$.
In \cite{EKLP_II}, there is only one weighting function, $w_0$, and the generating function is manipulated so that $w_0$ can be reused.  Such an approach would also apply to our setting, but it doesn't give us the results we want.

\begin{figure}
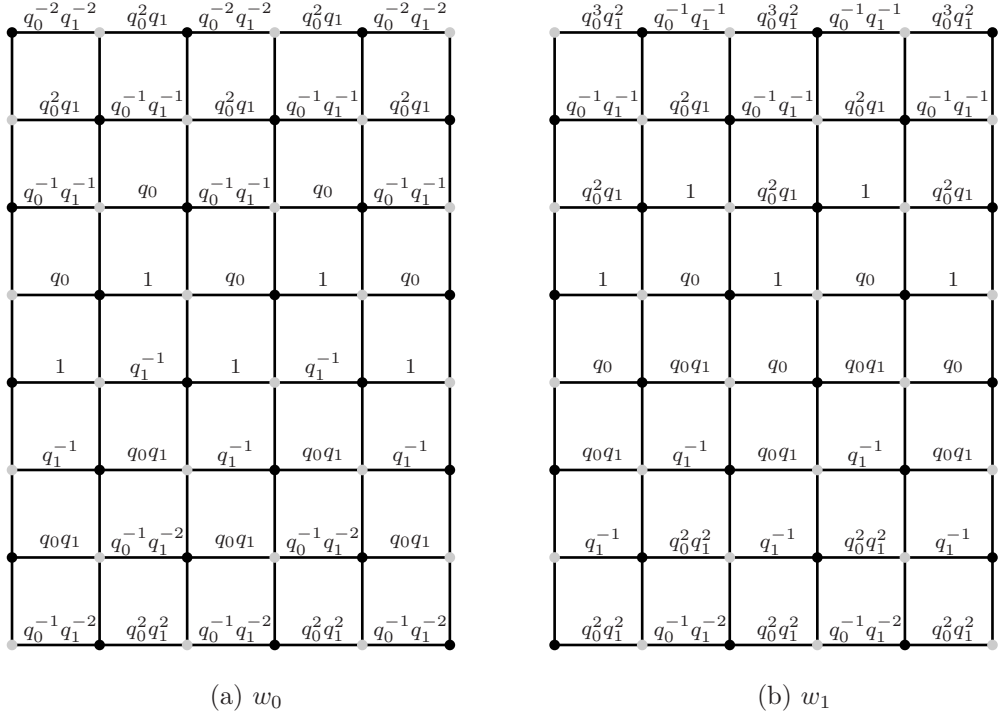

\caption{A comparison of the weightings $w_0$ and $w_1$}
\label{fig:shuffled_weights}
\rule{0in}{\baselineskip} \\
\subfloat[$w_0$]{ \input{unshuffled_weights.pstex_t} }
\subfloat[$w_1$]{\input{shuffled_weights.pstex_t} }
\end{figure}

\begin{lemma}
\label{lem:block_weight_lemma}
Let $d,d'$ be horizontal dimers, with $d'$ immediately north of $d$.  Then
\[
w_a(d)w_a(d') = \begin{cases}
q_0^{a+1}q_1^a & \text{if the block formed by }d, d' \text{ is odd},\\
q_0^{a}q_1^{a-1} & \text{if the block formed by }d, d' \text{ is even}.
\end{cases}
\]
\end{lemma}

\begin{proof}
When $a=0$, the lemma follows from the definition of $w_0$.  Now, suppose $a>0$.  If the block $d,d'$ is even, then $S^{-1}$ interchanges $d$ and $d'$, so  
\[
w_a(d)w_a(d') = w_{a-1}(S^{-1}(d))w_{a-1}(S^{-1}(d')) 
 = w_{a-1}(d')w_{a-1}(d) = q_0^{a}q_1^{a-1} \\
\]
by induction on $a$; otherwise, $(d, S^{-1}(d))$ and $(d', S^{-1}(d))$ are odd blocks under the alternate coloring, and we have
\begin{align*}
w_a(d)w_a(d') &= w_{a-1}(S^{-1}(d))w_{a-1}(S^{-1}(d')) \\
&= \frac{w_{a-1}(S^{-1}(d))w_{a-1}(d)w_{a-1}(d')w_{a-1}(S^{-1}(d'))}{w_{a-1}(d)w_{a-1}(d') } \\
&= \frac{(q_0^{a}q_1^{a-1})^2}{q_0^{a-1}q_1^{a-2}} = q_0^{a+1}q_1^a
\end{align*}
again by induction on $a$. 
\end{proof}

Next, we define what we mean by the weight of an odd-deficient or even-deficient dimer cover:

\begin{definition}
Let $\tilde \eta$ be an odd-deficient (or even-deficient) pyramid partition of length $n$.  Let $\pi$ be the pyramid partition obtained by filling in the missing odd (even) blocks of $\tilde \eta$ with pairs of vertical dimers.  Then we define
\[w_a(\tilde \eta) = w_a(\pi).\]
\end{definition}
If there are $m$ odd blocks in $\tilde \eta$, then
\begin{equation}
\label{eqn:weight_of_pi}
\sum_{\pi \text{ fills in } \tilde{\eta}}\negthickspace \negthickspace w_a(\pi) = (1+q_0^{a+1}q_1^a)^m w_a(\tilde{\eta})
\end{equation}
because each 
odd block of $\tilde{\eta}$ may be filled in two ways: we can use two vertical dimers (which each have weight 1) or we can use two horizontal dimers (which have a combined weight of $q_0^{a+1}q_1^{a}$ by Lemma~\ref{lem:block_weight_lemma}).  Similarly, if there are $m'$ odd blocks in $S(\tilde \eta)$, we have
\begin{equation}
\label{eqn:weight_of_S_pi}
\sum_{\pi' \text{ fills in } S(\tilde \eta)}\negthickspace \negthickspace w_{a+1}(\pi') = (1+q_0^{a+1}q_1^a)^{m'} w_{a+1}(S(\tilde \eta))
\end{equation}
As $\tilde \eta$ runs over $\tilde \mcP_n$, $S(\tilde \eta)$ runs over $\hat \mcP_{n+1}$. Also, Lemma~\ref{lem:odd_block_count} implies that $m-m'=n$, so combining Equations~(\ref{eqn:weight_of_pi}) and~(\ref{eqn:weight_of_S_pi}), we get
\begin{equation}
\label{eqn:general_shuffle_equation}
\sum_{\pi \in \mcP_n} w_a(\pi)
= (1+q_0^{a+1}q_1^{a})^n \sum_{\pi \in \mcP_{n+1}} w_{a+1}(\pi)
\end{equation}
Using Equation (\ref{eqn:general_shuffle_equation}) $k$ times, starting with $n=1$ and $a=0$, yields
\begin{equation}
\label{main_shuffle_result}
Z(1;q_0, q_1) = \left(\prod_{i=1}^k(1+q_0^iq_1^{i-1})^i\right)\sum_{\pi \in \mcP_{k+1}} w_{k}(\pi).
\end{equation}
As $k \rightarrow \infty$, the product on the right-hand side becomes $M(-q_1^{-1}, q_0q_1)^{-1}$, which is certainly good news, as this is one of the factors which appears in the statement of Theorem~\ref{thm:main}.  Next we need to try to understand the sum 
\[
\sum_{\pi \in \mcP_{k+1}} w_{k}(\pi)
\]
in the limit $k \rightarrow \infty$.

\section{Length-$\infty$ pyramid partitions}

\begin{figure}
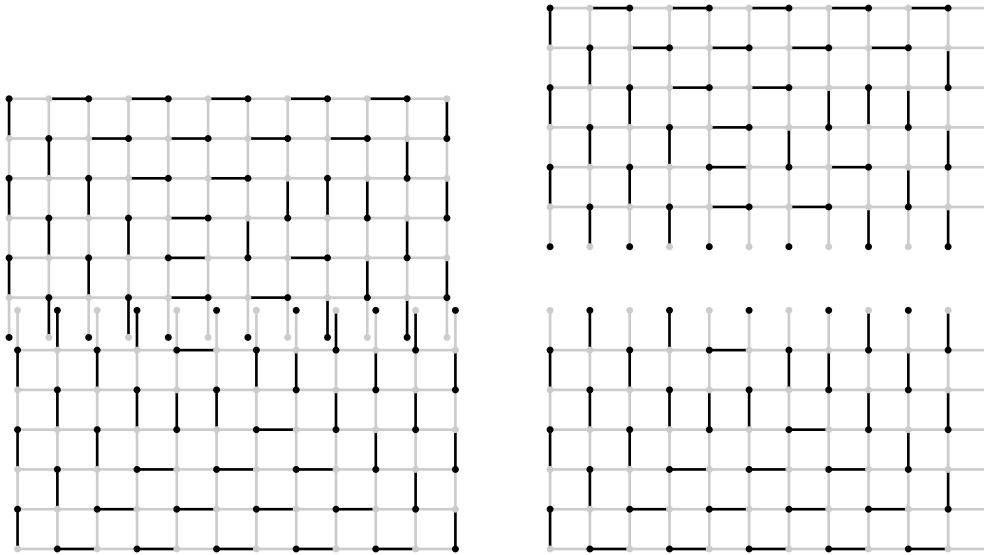

\caption{A pyramid partition, being split into two pieces which agree along the seam.}
\label{fig:splittem_up}
\rule{0in}{\baselineskip} \\
\input{splittem_up.pstex_t}
\input{really_splittem_up.pstex_t}
\end{figure}

In order to speak sensibly about the limit of the weighting functions $w_n$ as
$n$ gets large, we must shift our viewpoint slightly.  We shall split the
square lattice along the $x$ axis, giving us two half planes.  There are
infinitely many vertical edges which cross the $x$ axis; we shall include these
edges in both half-planes, and identify them.  A pyramid partition of length
$1$ therefore corresponds to two half-pyramid partitions which agree along the
``ragged'' edges of the two half-planes (see Figure~\ref{fig:splittem_up}).  Note that we don't quite have two matchings of the two graphs because the pendant edges (those that cross the $x$ axis) aren't necessarily in $\pi$.

This is a trivial change of viewpoint, but it allows us to shuffle the upper and lower half-planes independently.  When we are applying $S$ to the weights in the lower half-plane, let us imagine that we are travelling with the southbound weights.  From our new point of view, the northbound weights now move two units north, the ``westbound'' weights move northwest, and the ``eastbound'' weights move northeast.  Similarly, when we are applying $S$ to the upper half-plane, we are travelling with the northbound weights.

\begin{figure}
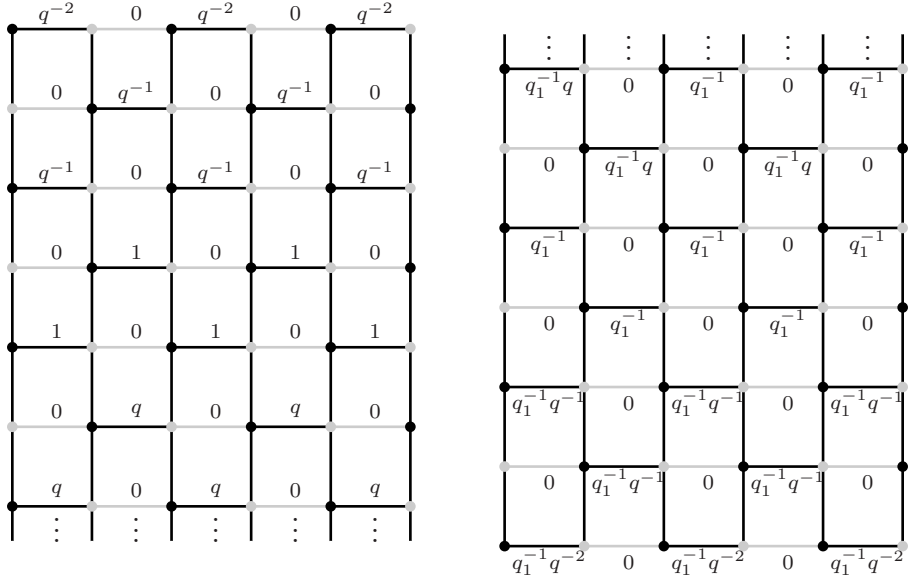

\caption{The weighting $w_\infty$, top and bottom pieces.}
\label{fig:w_infinity}
\rule{0in}{\baselineskip} \\
\input{w_infinity_top.pstex_t}
\input{w_infinity_bottom.pstex_t}
\end{figure}

Now it is clear what happens to the weight function $w_n$ as $n$ goes to infinity.  In the lower half-plane, nothing happens to the (now stationary) southbound edges at all.  However, the weights of the northbound edges get multiplied by $q_0q_1$.  
Let $q=q_0q_1$.  If we start with $n=1$ and shuffle $k$ times, the northbound edges are multiplied by $q^n$, so in the limit $n \rightarrow \infty$, they get weight zero.    In the same way, the southbound edges in the upper half plane get weight zero.  We call this weight function $w_{\infty}$; it is shown in Figure~\ref{fig:w_infinity}.

\begin{figure}
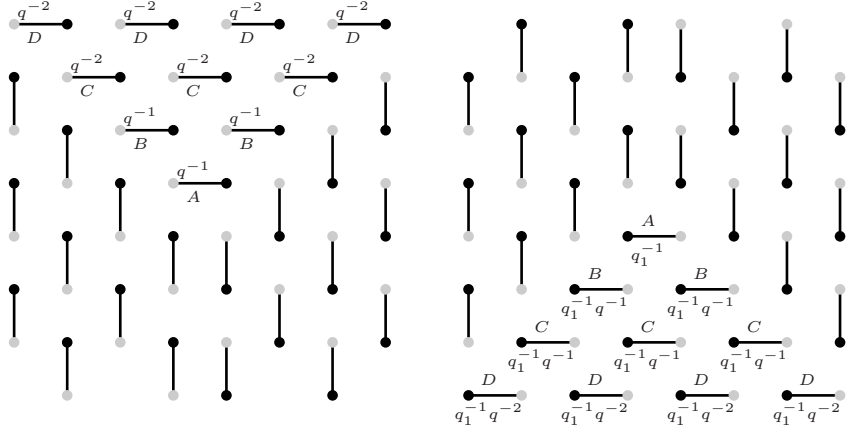

\caption{The empty room of length $\infty$, $\emptyroom_{\infty}$, top and bottom halves}
\label{fig:empty_infinity}
\rule{0in}{\baselineskip} \\
{\tiny
\input{empty_infinity_top.pstex_t}
\input{empty_infinity_bottom.pstex_t}
}
\end{figure}

We compute the weights of pyramid partitions in the same way as before: by normalizing by the weight of $\emptyroom_{\infty}$ (see Figure~\ref{fig:empty_infinity}). 
When we compute the sum $\sum_\pi w_{\infty}(\pi)$, we find that pyramid
partitions with southbound edges in the upper part, or northbound edges in the
lower part, get assigned weight zero.  Therefore, the only configurations $\pi$
that contribute to the sum $\sum w_\infty(\pi)$ are in fact perfect matchings
on the heavy edges in Figure~\ref{fig:w_infinity}, asymptotically identical to
the empty room of length infinity (see Figure~\ref{fig:empty_infinity}).  

Furthermore, if a dimer configuration of this type has horizontal edges arbitrarily far south in its upper half, or arbitrarily far north in its lower half, it also gets weight zero.  Thus the only dimer configurations that get nonzero weight under $w_{\infty}$ have a large frozen region of vertical dimers in the middle.  

\begin{definition}
A \emph{pyramid partition of length $\infty$} is a dimer configuration $\pi$ with $w_{\infty}(\pi) > 0$.
\end{definition}

  In order to determine whether pyramid partitions of length $\infty$ can be weight-enumerated in any sensible way, we should try to write down a set of elementary moves which can be applied to the empty room, sequentially, and are capable of generating all such $\pi$.   One such set is depicted in Figure~\ref{fig:infinity_moves}.  

One uses these elementary moves as follows.  Suppose we wish to construct a partition $\pi \in \mcP_{\infty}$.  Start with $\emptyroom_{\infty}$, and apply the ``infinite'' elementary move (a) until the frozen region in the middle is correct.  Then apply move (b) to the upper region and move (c) to the lower region until you have $\pi$.  

Note that move (a) deletes horizontal dimers from $\emptyroom_{\infty}$ symmetrically in pairs.  The first application of the move deletes the two dimers marked $A$ in Figure~\ref{fig:empty_infinity}; the next deletes two dimers marked $B$, and so on. 
Furthermore, the weight change of move (a) depends on where it is applied.  If two dimers marked $A$ are deleted, then the weight increases by $q_1q$; if two dimers marked $B$ are deleted, then the weight increases by $q_1q^2$, and so on.

\section{A weight-preserving bijection}
\label{sec:super_rigid}

We begin by defining super-rigid partitions, which are so named because they are a class of three-dimensional partitions whose generating function is the partition function for the Donaldson-Thomas theory of Calabi-Yau threefolds which come from super-rigid rational curves (see ~\cite{SR_DT}).

\begin{definition}
A \emph{Young diagram} is a finite subset of $(\bbZ_{\geq 0})^2$ which 
satisfies the following closure properties:
\begin{enumerate}
\item If $(x,y) \in \lambda$ and $x>0$, then $(x-1, y) \in \lambda$.
\item If $(x,y) \in \lambda$ and $y>0$, then $(x, y-1) \in \lambda$.
\end{enumerate}
\end{definition}

\begin{definition} Let $\lambda$ be a Young diagram. The \emph{leg of shape $\lambda$} is the set 
\[
L_\lambda = \lambda \times \bbZ_{\geq 0} \subseteq (\bbZ_{\geq 0})^3.
\]
A \emph{three-dimensional partition asymptotic to $\lambda$} is a set $\pi$ satisfying $L_\lambda \subseteq \pi \subseteq (\bbZ_{\geq 0})^3$ satisfying the following properties: 
\begin{enumerate}
\item The set $\pi \setminus L_\lambda$ is finite.
\item If $(x,y,z) \in \pi$ and $x>0$, then $(x-1, y, z) \in \pi$.
\item If $(x,y,z) \in \pi$ and $y>0$, then $(x, y-1, z) \in \pi$.
\item If $(x,y,z) \in \pi$ and $y>0$, then $(x, y, z-1) \in \pi$.
\end{enumerate}
We also define the \emph{size of $\pi$}, written $|\pi|$, to be the cardinality of the set $\pi \setminus L_{\lambda}$.
\end{definition}

If $\pi$ is a three dimensional partition asymptotic to $\lambda$, we informally call the elements of $\pi$ ``boxes''; one can think of $\pi$ as a stack of boxes in the corner of a large room which has one ``baseboard'' whose cross-section is $\lambda$.

\begin{definition}
A \emph{super-rigid partition} is a triple $(\pi_0, \lambda, \pi_{\infty})$, where $\pi_0$ and $\pi_\infty$ are three-dimensional partitions asymptotic to $\lambda$.  
\end{definition}

\begin{lemma}[Lemma 2.9 of~\cite{SR_DT}]
\label{DT_lemma}
  Give the super-rigid partition $(\pi_0, \lambda, \pi_{\infty})$ the weight $z^{|\lambda|}q^N$, where 
\[
N = |\pi_0| + |\pi_{\infty}| + \sum_{i,j \in \lambda} (i+j+1).
\]
The generating function for super-rigid partitions under this weighting scheme is 
\[
Z_X(z,q) = M(1,q)^2 M(-z, q)^{-1}.
\]
\end{lemma}
\begin{definition} $Z(\infty; q_0, q_1) = Z_X(q_1, q_0q_1).$
\end{definition}
\begin{figure}
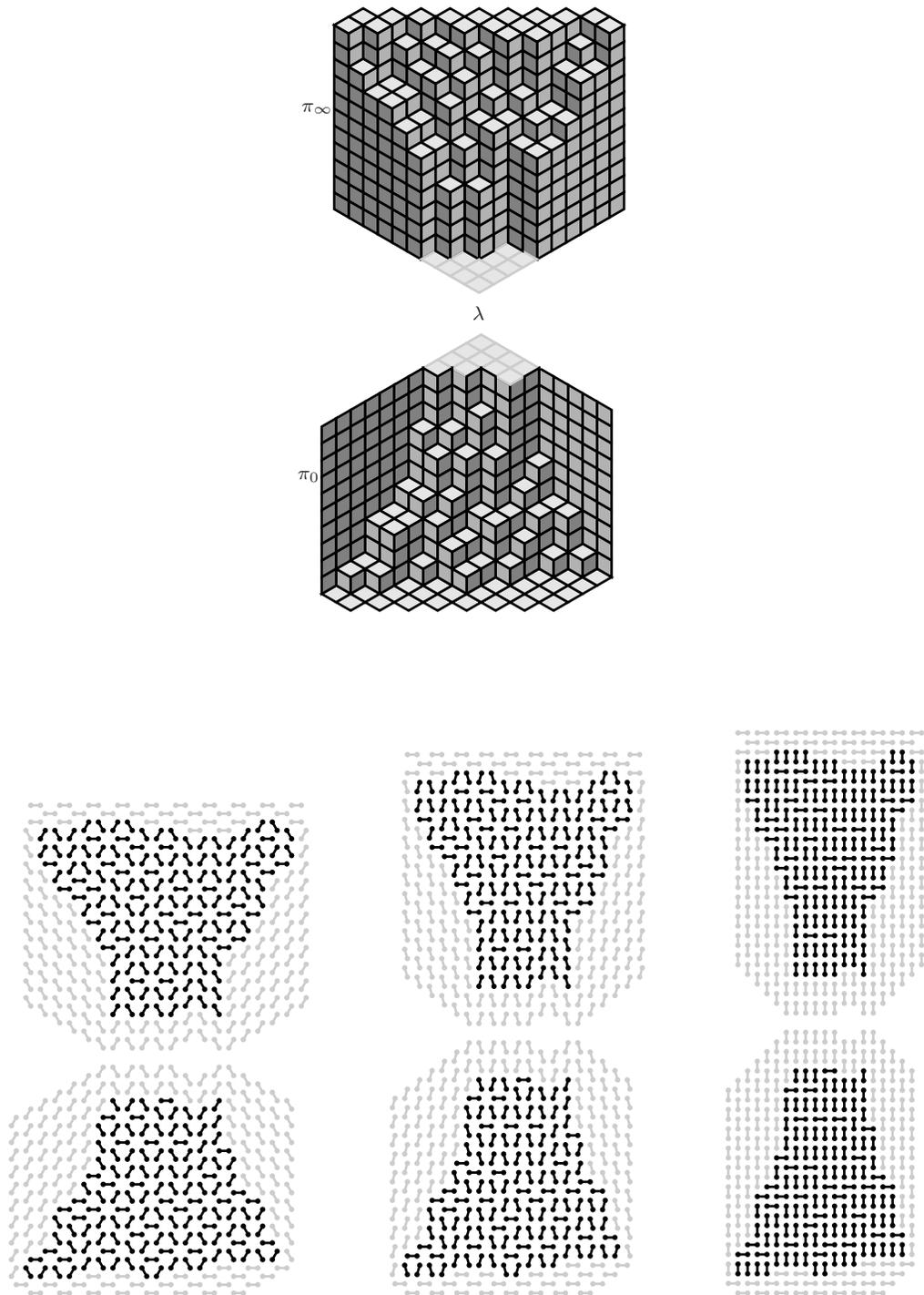

\caption{A super-rigid partition $(\pi_0, \lambda, \pi_{\infty})$ becoming a pyramid partition of length $\infty$}
\label{fig:super_rigid}
\raisebox{1.5in}{
\begin{tabular}{c}
\thinspace \thinspace \thickspace \thickspace 
\input{super_rigid_shade_top.pstex_t} \\
\input{super_rigid_shade_bottom.pstex_t} 
\end{tabular}
}  \\ \rule{0in}{2 \baselineskip} \\ 
\input{super_rigid_hex.pstex_t}
\input{super_rigid_nearpyramid.pstex_t}
\input{super_rigid_pyramid.pstex_t}
\end{figure}

There is a ``folklore'' correspondence between 3D partitions and dimer covers
of the hexagon lattice: if we view a 3D partition from far away along the line
$x=y=z$, it appears to be a tiling of the plane by lozenges.  Replacing each of
these lozenges with a dimer, we get a dimer cover of the hexagon lattice.  A
simple reorientation of the edges of the hexagon lattice shows that it is the
same as the ``brickwork'' lattices defined by the heavy lines of
Figure~\ref{fig:w_infinity}.  

Let us apply this observation to create a correspondence between super-rigid
partitions and pyramid partitions of length $\infty$.  Starting with $(\pi_0,
\lambda, \pi_{\infty})$, replace both $\pi_0$ and $\pi_{\infty}$ by their dimer
versions, and then reorient all of the edges so that the dimers fit onto the
brickwork lattice (see Figure~\ref{fig:super_rigid}).  The fact that $\pi_0$
and $\pi_{\infty}$ share a common asymptotic leg $\lambda$ causes the frozen
region of vertical dimers to appear in the middle of the figure.  The
correspondence is clearly bijective, and with a little care, we can make this
bijection weight-preserving.

Consider any super-rigid partition $(\pi_0, \lambda, \pi_{\infty})$.  We can construct this partition from the empty super-rigid partition $(\emptyset, \emptyset, \emptyset)$ using the following three elementary moves:

\begin{itemize} 
\item[(a)] Add a log of boxes in position $(i,j)$ to the asymptotic leg, with weight $q_1q^{i+j+1}$.  Repeat until we have constructed $\lambda$.
\item[(b)] Add a box to the left end of the partition, with weight $q$.  Repeat until we have constructed the super-rigid partition $(\pi_0, \lambda, \emptyset)$.
\item[(c)] Add a box to the right end of the partition, with weight $q$.  Repeat until we have constructed $(\pi_0, \lambda, \pi_{\infty})$.
\end{itemize} 

\begin{figure}
\caption{Elementary moves for generating elements of $\mcP_{\infty}$ from $\emptyroom_\infty$}
\label{fig:infinity_moves}
\rule{0in}{\baselineskip} \\
\subfloat[Middle region]{\input{infinity_moves.pstex_t}}
\subfloat[Upper region]{\input{infinity_move_upper.pstex_t}}
\subfloat[Lower region]{\input{infinity_move_lower.pstex_t}}
\end{figure}

A partition constructed in this manner will be weighted correctly to contribute to $Z(\infty; q_0, q_1)$.
Note that we have deliberately chosen these moves to have the same names as those in Figure~\ref{fig:infinity_moves}.  
Define a bijection \[
\Phi:\mcP_{\infty} \rightarrow \{\text{super-rigid partitions}\}
\]
as follows:
given $\pi \in \mcP_{\infty}$, determine a set of elementary moves to construct
$\pi$ from $\emptyroom_{\infty}$, and then use the corresponding moves in the same
order to create a super-rigid partition.  This super-rigid partition is
$\Phi(\pi)$. 

Since each of these elementary moves affects the weight in the same manner as the corresponding move on pyramid partitions, $\Phi$ is weight-preserving.  Thus $\Phi$ also preserves the generating functions:
\[
\sum_{\pi \in \mcP_{\infty}} w_{\infty}(\pi) = Z(\infty; q_0, q_1).
\]
In the limit $n \rightarrow \infty$, Equation~\ref{main_shuffle_result} now says
\[
Z(1;q_0, q_1) = \left(\prod_{i=1}^{\infty}(1+q_0^iq_1^{i-1})^i\right) Z(\infty; q_0, q_1)
\]
which proves Theorem~\ref{thm:main}.

\section{The generating function for general $n$}

Next, we shall use the same argument to calculate $Z(n; q_0, q_1)$.  Applying Equation~(\ref{eqn:general_shuffle_equation}) $k$ times, starting at $a=0$ but leaving $n$ arbitrary, we get
\begin{equation}
\label{general_shuffle_result}
Z(n;q_0, q_1) = \left(\prod_{i=1}^k(1+q_0^iq_1^{i-1})^{i+n-1}\right)\sum_{\pi \in \mcP_{k+n}} w_{k}(\pi).
\end{equation}
Taking the limit as $k$ approaches infinity, we again get a sum over pyramid partitions of length $\infty$, but with a slightly modified weight function $w_{\infty}^{n}$:
\begin{equation}
\label{eqn:intermediate_weight_n_step}
Z_n(;q_0, q_1) = \left(\prod_{i=1}^{\infty}(1+q_0^iq_1^{i-1})^{i+n-1}\right)\sum_{\pi \in \mcP_{\infty}} w_{\infty}^n(\pi),
\end{equation}
$w_{\infty}^{n}$ has the property that the elementary move of type (a) carries the weight $q_1q^{i+j+n}$.  This means that the corresponding super-rigid partition $(\pi_0, \lambda, \pi_{\infty})$ has weight $q_1^{\lambda}q^{N(n)}$, where
\[
N(n) = |\pi_0| + |\pi_{\infty}| + (n-1)|\lambda| + \sum_{i,j \in \lambda} (i+j+1).
\]
We only need a slight modification to the argument of~\cite{SR_DT} to compute the sum on the right-hand side of Equation~\ref{eqn:intermediate_weight_n_step}.  We begin with the one-leg formula for the topological vertex (see \cite{CRYSTALS}) , which states that
\[
\sum_{\pi \text{ asymp. to } \lambda} q^{|\pi|}
= M(q) q^{\binom{\lambda}{2}}s_{\lambda^t}(q),
\]
where $\binom{\lambda}{2} = \sum_{\lambda_i \in \lambda} \binom{\lambda_i}{d}$, $\lambda^t$ denotes the transpose of $\lambda$, and $s_{\lambda^t}(q)$ denotes the principal specialization of the Schur function.  We have
\begin{align*}
\sum_{\pi \in \mcP_{\infty}} w_{\infty}^n(\pi)
&= \sum_{\lambda} \sum_{\pi_0, \pi_{\infty} \rightarrow \lambda} q^{N(n)}q_1^{\lambda} \\
&= M(1,q)^2 \sum_{\lambda}q_1^{\lambda}q^{(n-1)\lambda + 
\binom{\lambda}{2} + \binom{\lambda^t}{2} + \sum_{(i,j) \in \lambda} i+j-1}
s_{\lambda^t}(q)s_{\lambda}(q) \\
&= M(1,q)^2 \sum_{\lambda}q_1^{\lambda}q^{n|\lambda|}
s_{\lambda^t}(q)s_{\lambda}(q) \\
&= M(1,q)^2 \prod_{i,j = 1}^{\infty} (1+q_1 q^{i+j+n-2}) \\
&= M(1,q)^2 \prod_{k = 1}^{\infty} (1 + q_0^kq_1^{k+1})^{\max(k-n+1, 0)}.
\end{align*}
This proves Theorem~\ref{thm:generalization_of_main}. $\hfill \square$

\section{Future Work}

There are several possible lines of research suggested by the techniques and results of this paper:

\begin{enumerate}
\item  The shuffling procedure still works for certain pyramid partitions which are not asymptotic to $\emptyroom_n$.  In particular, we can allow pyramid partitions to have up to four asymptotic legs, pointing NW, NE, SW, and SE, whose shapes are given by partitions $\lambda_{NW}, \lambda_{NE}, \lambda_{SW}, \lambda_{SE}$.  It seems possible that we could compute the generating function for such configurations using the full topological vertex formula of~\cite{CRYSTALS}.  Such a result might shed some light on flop transitions in topological string theory.

\item It may be possible to compute a somewhat more refined generating function,using $2n$ variables rather than just two.  This would have the effect of introducing diagonal ``stripes'' on the alternate layers of the pyramid partition.  Such a count is done in~\cite{ORBIFOLD_DT} using vertex operator methods.

\item This paper shows that there is a direct link between the Donaldson-Thomas partition function of the conifold, $Z(1, q_0, q_1)$, and the Donaldson-Thomas partition function of the resolution, $Z(\infty, q_0, q_1)$.  We have proven~\cite{ORBIFOLD_DT} that there is a similar relationship between the Donaldson-Thomas partition function of the orbifold $\bbC^3/G$, (where $G$ is a finite Abelian subgroup of $SO(3)$), and the Donaldson-Thomas partition function of its resolution; unfortunately, the methods of~\cite{ORBIFOLD_DT} do not suggest why this should be.  We can attempt to look for this type of relationship between other singular threefolds and their resolutions.

\end{enumerate}

\section*{Acknowledgements}
The author would like to thank Dr. Bal\'azs Szendr\H{o}i for posing this problem; Dr. Richard Kenyon for a fruitful conversation which led to its solution; and Dr. Jim Bryan for his careful proofreading and helpful suggestions.  

\bibliographystyle{plain}
\bibliography{conifold}

\end{document}